\documentclass[12pt]{article}

\usepackage[english]{babel}
\usepackage{amssymb,amsmath}

\newtheorem{theorem}{Theorem }[section]

\newtheorem{definition}[theorem]{Definition}

\newcommand{\qed}{\hspace*{\fill}$\Box$}

\title{Projective Two-Weight Sets of Denniston Type} 

\author{Stefaan De Winter}

\date{The National Science Foundation}
\begin{document}

\maketitle

\abstract{We construct two-weight sets in PG$(3n-1,q)$, $n\geq2$ with the same weights as those that would arise from the blow-up of a maximal $q$-arc in PG$(2,q^n)$. The construction is of particular interest when $q$ is odd, as it is well known that no maximal arcs in PG$(2,q^n)$ exist in that case.}

	\section{Introduction}

	{\it Projective two-weight sets of type $(a,b)$}, that is, point sets in PG$(n,q)$ such that each hyperplane intersects the point set in either $a$ or $b$ points, play a central role in finite geometry due to their equivalence with linear two-weight codes and strongly regular Cayley graphs on an elementary Abelian group. These equivalences, as well as an early survey of projective two-weight sets, are described in great detail in the seminal paper ``The Geometry of Two-Weight Codes'' by Calderbank and Kantor \cite{CK}. Among projective two-weight sets an important role is reserved for the {\it maximal arcs} in PG$(2,q)$ which apart from the above mentioned connections also give rise to so-called partial geometries. We provide a formal definition:
	
	\begin{definition}
	A maximal $d$-arc $\mathcal{K}$ in PG$(2,q)$ is a nonempty set of points in PG$(2,q)$ with the property that every line intersects $\mathcal{K}$ in either $0$ or $d$ points with $1<d<q$.
	\end{definition}
	
	It follows immediately from the definition that necessarily $d\mid q$. When $q=2^n$ Denniston \cite{Denniston} showed that this condition is sufficient for the existence of a maximal $d$-arc, that is, for every $0<k<n$ there exists a maximal $2^k$-arc in PG$(2,2^n)$. However, a by now classical result by Ball, Blokhuis and Mazzocca \cite{BBM} states that no maximal arc in PG$(2,q)$ can exist whenever $q$ is odd. Maximal arcs have been studied extensively from both a constructive perspective (e.g. \cite{Thas}, \cite{Mathon} \dots), a coding theoretical perspective (e.g. \cite{Ding},\cite{Tonchev},\dots), an algebraic perspective (e.g. \cite{Qing1}, \cite{pent},\dots),  as well as an enumerative perspective (e.g. \cite{ball}, \cite{Maes},\dots).
	
	Consider a maximal $d$-arc in PG$(2,q^n)$. As the vector space $V(3,q^n)$ underlying this plane can be interpreted as a vector space $V(3n,q)$, we see that any point set in PG$(2,q^n)$ naturally gives rise to a set of disjoint $(n-1)$-dimensional subspaces of PG$(3n-1,q)$ (the so-called ``blow-up'' of PG$(2,q)$). It is easy to see that the blow-up of a maximal $d$-arc in PG$(2,q)$ is a two-weight set with weights $((q^nd-q^n+d)\frac{q^{n-1}-1}{q-1},d\frac{q^{n}-1}{q-1}+(q^nd-q^n)\frac{q^{n-1}-1}{q-1})$ in PG$(3n-1,q)$. It is now natural to wonder if two-weight sets with such weights exist when $q$ is odd (for $n\geq2$). Note that there is no requirement though that such set would consist of the union of disjoint $(n-1)$-spaces. The example the author was aware of for a long time and that formed the motivation for this research, is a projective two-weight set with weights $(21, 30)$ in PG$(5,3)$ constructed by computer by Mathon and put in a theoretical framework (from a projective geometric perspective) in \cite{perp} by De Clerck, Delanote, Hamilton and Mathon. The parameters of this two-weight set correspond to the blow-up of a non-existing maximal $3$-arc in PG$(2,9)$. What makes this example especially unique and interesting is that it is actually the union of $21$ disjoint lines.  The fact that two-weight sets can be studied from both a geometric, coding theoretic and graph theoretic perspective complicates the literature search for other examples as these communities do not always use the same ``language'' or communicate their results to one another. However, when looking for small examples the table of small strongly regular graphs by Brouwer \cite{Brouwer} provides an access point for looking for possible other examples. This way the author discovered that two-weight sets with weights $(21, 30)$ in PG$(5,3)$ were known before in a coding theoretic setting, constructed both by Gulliver in \cite{Gulliver} and Kohnert in \cite{Kohnert}. During the preparation of this manuscript these were the only examples the author was aware of. However, shortly before submission of this article the author came across a paper once given to him (and unfortunately not read at that time) by Jurgen Bierbrauer. In this paper \cite{Bierbrauer} Bierbrauer and Edel, through an intricate coding theoretic approach, construct a family of two-weight codes  that correspond to two-weight sets with the weights we are interested in. Their construction actually yields such two-weight sets for all allowable $d$ and $q$. However, their construction provides no geometric insight into the structure of these two-weight sets but merely establishes their existence through the equivalence with two-weight codes. Also, it is very unfortunate that this paper seems to be largely unknown in either the coding theory, graph theory or finite geometry community, having received only three citations as per MathSciNet. We hope that the current paper can help in making the Bierbrauer-Edel paper more broadly known. \medskip
\medskip

	In this article we will provide an elegant geometric construction of two-weight sets in PG$(3n-1,q)$ with weights $\frac{q^{2n}-q^{2n-1}+q^n-q}{q-1}$ and $\frac{q^{2n}-q^{2n-1}-q^{n+1}+2q^n-q}{q-1}$, that is, the weights corresponding to the blow-up of a maximal $q$-arc in PG$(2,q^n)$. We will also provide a more algebraic description of these sets. We will end the paper with some open questions.

	\section{The Cossidente-Storme Observation}
	
	In this short section we describe an observation originally made by Cossidente and Storme in \cite{CS} which is essential to our construction. Early versions of their result appear in Ebert \cite{Ebert} and Kestenband \cite{Kestenband}. Let the projective space PG$(2n-1,q)$ be represented by $\mathbb{F}_{q^{2n}}$ $\pmod{ \mathbb{F}_q}$. This means that the points of PG$(2n-1,q)$ are the non-zero elements of $\mathbb{F}_{q^{2n}}$ and two elements $x$ and $y$ represent the same point if and only if $x/y \in \mathbb{F}_q$. Now let $\beta$ be a primitive element of $\mathbb{F}_{q^{2n}}$. Then the mapping $\sigma: x\mapsto \beta x$ with $x\in\mathbb{F}_{q^{2n}}^*$ generates a Singer cycle of PG$(2n-1,q)$, that is, it has a unique orbit on the points of PG$(2n-1,q)$. We can now state the necessary results from \cite{CS}.
	
	\begin{theorem}
	For every $a\in\mathbb{F}_{q^{n}}^*$  the set $\mathcal{Q}_a=\{x\in \mathrm{PG}(2n-1,q) \mid \mathrm{Tr}(ax^{q^n+1})=0\}$, with Tr the usual trace map from $\mathbb{F}_{q^{n}}$ to $\mathbb{F}_q$, is an elliptic quadric of PG$(2n-1,q)$.
	\end{theorem}
	
	This is Theorem 3.2 in \cite{CS}. It is easily seen that $\mathcal{Q}_a=\mathcal{Q}_b$ if and only if $a/b\in\mathbb{F}_q$. Hence one obtains $(q^n-1)/(q-1)$ distinct elliptic quadrics and can observe the following:
	
	\begin{theorem}
	The elliptic quadrics $\mathcal{Q}_a$ form, through the parameter $a$, a projective space $\mathrm{PG}$$(n-1,q)$. Any linear combination $\lambda\mathcal{Q}_a+\mu\mathcal{Q}_b$, $\lambda,\mu\in\mathbb{F}_q$, of two elliptic quadrics $\mathcal{Q}_a$ and $\mathcal{Q}_b$ defines a new elliptic quadric $\mathcal{Q}_{\lambda a+\mu b}$.
	\end{theorem}
	
	This is Remark 3.5 in \cite{CS} and is essential for our construction. Finally we need some special orbits under a subgroup of our above described Singer cycle.
	
	\begin{theorem}\label{key}
	Let $\xi=\beta^{(q^n-1)/(q-1)}$, with $\beta$ a primitive element of $\mathbb{F}_{q^{2n}}$. Then an orbit in $\mathrm{PG}(2n-1,q)$ under $<\xi>$ has size $q^n+1$ and is either contained in or disjoint from $\mathcal{Q}_a$, $a\in\mathbb{F}_{q^{n}}$. Furthermore, each such orbit is either a cap (when $n$ is even) or a union of disjoint lines (when $n\geq3$ is odd) that is the intersection of $n-1$ linearly independent $\mathcal{Q}_a$.
	\end{theorem}
	
	This is theorem 3.1 combined with theorems 3.7 and 3.8 in \cite{CS}. As a result we can view the elliptic quadrics $\mathcal{Q}_a$ as the points of a PG$(n-1,q)$ (see above) and the orbits under $<\xi>$ as the hyperplanes of this projective space.
	
	\section{The Two-Weight Set $\mathcal{C}$}
	
	Let $n\geq2$. Let $\Lambda:=$PG$(2n-1,q)$ be embedded in PG$(3n-1,q)$. Let $\Pi:=$PG$(n-1,q)$ be embedded in PG$(3n-1,q)$ disjoint from $\Lambda$.
	
	Let $\mathcal{Q}_i$, $i=1,\dots, (q^n-1)/(q-1)$, be a set of non-singular elliptic quadrics covering $\Lambda$ as described in the previous section, such that the $\mathcal{Q}_i$ form the points of a PG$(n-1,q)$, say $\Gamma$.
	
	Let $\delta$ be an anti-isomorphism between $\Gamma$ and $\Pi$. Assume that the hyperplanes $\pi_i$ of $\Pi$ are labeled in such a way that $\delta(\mathcal{Q}_i)=\pi_i$.
	
	By Theorem \ref{key} we know that the intersection of $n-1$ linearly independent quadrics $\mathcal{Q}_i$ contains exactly $q^n+1$ points (either a cap or union of disjoint lines depending on whether $n$ is even or odd). Label these intersections as $I_i$ for $i=1,\dots, \frac{q^n-1}{q-1}$ (note that this is the correct number of these intersections as these correspond to the hyperplanes of $\Gamma$). Hence every point $p$ of $\Pi$ is associated under $\delta$ in a unique way with one of the $I_i$. Assume the points of $\Pi$ are labeled $p_i$ in such a way that $p_i$ corresponds to $I_i$ in the above described correspondence.

	Define the baseless cone $C_i$ as $C_i:=p_iI_i\setminus\{I_i\}$. Note that all these cones are mutually disjoint. This follows from the fact that the $I_i$ are disjoint, which in turns follows from their construction through a subgroup of a Singer group.
		
	\medskip
	
	{\bf Theorem.}\label{main} {\it The union of cones $\mathcal{C}=\bigcup_i C_i$ forms a two-weight set in {\rm PG}$(3n-1,q)$, $n\geq2$, with weights $\frac{q^{2n}-q^{2n-1}+q^n-q}{q-1}$ and $\frac{q^{2n}-q^{2n-1}-q^{n+1}+2q^n-q}{q-1}$.}
	\medskip
	
	{\bf Proof.} There are three types of hyperplanes to consider. Those containing $\Lambda$, those containing $\Pi$, and those intersecting both $\Lambda$ and $\Pi$ in a hyperplane. Before discussing these three cases we observe that obviously there are $\frac{q^n-1}{q-1}$ quadrics $\mathcal{Q}_i$, each of these $\mathcal{Q}_i$ contains $\frac{(q^n+1)(q^{n-1}-1)}{q-1}$ points, and each point of $\Lambda$ is contained in $\frac{q^{n-1}-1}{q-1}$ of the quadrics $\mathcal{Q}_i$.
	\begin{itemize}
		\item Let $\alpha$ be a hyperplane containing $\Lambda$. Then $\alpha$ intersects $\Pi$ in some hyperplane $\pi$, and $\alpha\cap\mathcal{C}$ consists exactly of the cones $C_i$ for which $p_i\in \pi$. It follows that \begin{eqnarray*}|\alpha\cap\mathcal{C}|&=&\frac{q^{n-1}-1}{q-1}(q-1)(q^n+1)+\frac{q^{n-1}-1}{q-1} \\ &=& \frac{q^{2n}-q^{2n-1}-q^{n+1}+2q^n-q}{q-1}.\end{eqnarray*}
			
		\item Let $\alpha$ be a hyperplane containing $\Pi$. Then $\alpha$ intersects $\Lambda$ in a hyperplane $\lambda$. Now each point of $\lambda$ is contained in unique $I_i$ and hence belongs to a unque line in $p_iI_i$. It follows that \begin{eqnarray*}|\alpha\cap\mathcal{C}|&=&\frac{q^{2n-1}-1}{q-1}(q-1)+\frac{q^n-1}{q-1} \\ &=& \frac{q^{2n}-q^{2n-1}+q^n-q}{q-1}.\end{eqnarray*}
		
		\item Finally we consider the case where $\alpha$ is a hyperplane that intersects both $\Lambda$ in a hyperplane $\lambda$ and $\Pi$ in a hyperplane $\pi_i$. There are two case to consider depending on how $\mathcal{Q}_i$ intersects $\lambda$. Recall that $\delta(\mathcal{Q}_i)=\pi_i$. We observe that $\lambda\cap\mathcal{Q}_i$ can either be a parabolic quadric $Q(2n-2,q)$ or a cone $pQ^-(2n-3,q)$. In either case, while we do not know how a specific $I_j$ intersects $\lambda$ we can still derive $|<\lambda,\pi_i>\cap\ \mathcal{C}|$. Every point $p$ in the intersection of $Q_i$ and $\lambda$ is joined in $\mathcal{C}$ by a unique line to some point of $\pi_i$. Note that the latter follows from the fact that $\delta$ is an anti-automorphism, implying that the image under $\delta^{-1}$ of a point of $\pi_i$ must be a ``hyperplane'' through $\mathcal{Q}_i$, that is, some $I_j$ contained in $\mathcal{Q}_i$.  Furthermore, no other point of $\lambda$ can belong to any of the cones $p_jI_j$ for $p_j\in\pi_i$. Hence $|<\lambda,\pi_i>\cap\ \mathcal{C}|=\frac{q^{2n-2}-1}{q-1}(q-1)+\frac{q^{n-1}-1}{q-1}$ if $\lambda\cap\mathcal{Q}_i$ is a parabolic quadric $Q(2n-2,q)$, and $|<\lambda,\pi_i>\cap\ \mathcal{C}|=(1+q\frac{(q^{n-1}+1)(q^{n-2}-1)}{q-1})(q-1)+\frac{q^{n-1}-1}{q-1}$ if $\lambda\cap\mathcal{Q}_i$ is a cone $pQ^-(2n-3,q)$. Now the hyperplane $\alpha$ contains $<\lambda,\pi_i>$ and intersects each of the lines $<r,u>$, $r\in\Lambda\setminus\lambda$, $u\in\Pi\setminus\pi_i$ in a unique point. Furthermore every point of $\alpha\cap\mathcal{C}$ not in $<\lambda,\pi_i>$ must be contained in such line $<r,u>$. Also recall that every point $r\in\Lambda$ is joined to a unique point of $\Pi$ through a line of one of the cones $C_j$. Start with the case where $\lambda\cap\mathcal{Q}_i$ is a parabolic quadric $Q(2n-2,q)$. The points of $\pi_i$ are joined to a total of $\frac{q^{n-1}-1}{q-1}(q^n+1)$ points of $\Lambda$ through a line of one of the cones $C_j$. As noted above $\frac{q^{2n-2}-1}{q-1}$ of these points belong to $\lambda$. Hence exactly  $\frac{q^{n-1}-1}{q-1}(q^n+1)-\frac{q^{2n-2}-1}{q-1}$ of the points of $\Lambda\setminus\lambda$ are not on one of the lines  $<r,u>$, $r\in\Lambda\setminus\lambda$, $u\in\Pi\setminus\pi_i$ belonging to a cone $C_j$, while all other points of $\Lambda\setminus\lambda$ are on a unique such line. As $\alpha$ intersects each of these lines in a unique point off of $\Lambda$ and $\Pi$ we obtain \begin{eqnarray*}|\alpha\cap\mathcal{C}|&= & \left(\frac{q^{2n-2}-1}{q-1}(q-1)+\frac{q^{n-1}-1}{q-1}\right) \\ & & +\left(q^{2n-1}-(\frac{q^{n-1}-1}{q-1}(q^n+1)-\frac{q^{2n-2}-1}{q-1})\right) \\ &=& \frac{q^{2n}-q^{2n-1}+q^n-q}{q-1}.\end{eqnarray*}
		Finally consider the case where $\lambda\cap\mathcal{Q}_i$ is a cone $pQ^-(2n-3,q)$. A similar argument as above, replacing the number of points $\frac{q^{2n-2}-1}{q-1}$ of a parabolic quadric $Q(2n-2,q)$ by the number of points $1+q\frac{(q^{n-1}+1)(q^{n-2}-1)}{q-1}$ of a cone  $pQ^-(2n-3,q)$, yields \begin{eqnarray*}|\alpha\cap\mathcal{C}|&= & \left((1+q\frac{(q^{n-1}+1)(q^{n-2}-1)}{q-1})(q-1)+\frac{q^{n-1}-1}{q-1}\right) \\ & & +\left(q^{2n-1}-\left(\frac{q^{n-1}-1}{q-1}(q^n+1)-(1+q\frac{(q^{n-1}+1)(q^{n-2}-1)}{q-1})\right)\right) \\ &=& \frac{q^{2n}-q^{2n-1}-q^{n+1}+2q^n-q}{q-1}.\end{eqnarray*}
		
	\end{itemize}
	This concludes the proof. \qed
	\bigskip
	
	Note that the weights of the above constructed two-weight sets correspond to those that would arise if one were to blow up a maximal $q$-arc in PG$(2,q^n)$ to PG$(3n-1,q)$. Of course for odd $q$ no such maximal arcs exist, yet the above construction works for all odd $q$ as well as for even $q$. 
	\medskip

	Note that in the case $n=2$ we do not need an anti-isomorphism between the so-called elliptic fibration of PG$(3,q)$ and the line PG$(1,q)$. Any bijection will let the construction work. This is no longer the case for $n\geq3$.
	\medskip
	
	\section{An Algebraic Construction}
	
	In this section we show that the set $\mathcal{C}$ admits an interesting group of automorphisms from which it can be reconstructed. We first coordinatize PG$(3n-1,q)$ in a somewhat uncommon way. Let the points of  PG$(3n-1,q)$ be coordinatized by 2-tuples $(x,y)\neq (0,0)$, with $x\in \mathbb{F}_{q^{2n}}$ and $y\in \mathbb{F}_{q^n}$, determined up to scalar multiple in $\mathbb{F}_q$; and where the tuples $(x,0)$ coordinatize $\Lambda$, and the tuples $(0,y)$ coordinatize $\Pi$ in the obvious way. 
		
	Let $\beta$ be a primitive element of $\mathbb{F}_{q^{2n}}$ and let $\gamma$ be a primitive element of $\mathbb{F}_{q^n}$. We further can assume that $\mathbb{F}_{q^n}\leq\mathbb{F}_{q^{2n}}$ and that $\beta$ and $\gamma$ are chosen such that $\beta^{q^n+1}=\gamma$. Now consider the group 
	$$G=<\left(\begin{array}{cc}\beta & 0 \\ 0 & \gamma\end{array}\right)>.$$
	
	This group naturally acts as a group of automorphisms of PG$(3n-1,q)$ through left multiplication with the points $(x,y)$ viewed as column vectors. 
	
	The points of $\Pi$ are represented by the 2-tuples $$(0,1),(0,\gamma),(0,\gamma^2),\dots,(0,\gamma^{\frac{q^n-1}{q-1}-1}),$$
	
	whereas the points of $\Lambda$ are represented by $$(1,0),(\beta,0),(\beta^2,0),\dots,(\beta^{\frac{q^{2n}-1}{q-1}-1},0).$$
	
	The stabilizer in $G$ of a point of $\Pi$ is easily seen to be given by the group $$S=<\left(\begin{array}{cc}\beta^{\frac{q^n-1}{q-1}} & 0 \\ 0 & \gamma^{\frac{q^n-1}{q-1}}\end{array}\right)>.$$

	The orbits in $\Lambda$ of $S$ are the sets $\{ (\beta^i\beta^{j\frac{q^n-1}{q-1}},0) \mid j=0,1,\dots,q^n \}$ for $i=0,1,\dots,\frac{q^n-1}{q-1}-1$. By Theorem \ref{key} we know these orbits are exactly the sets $I_i$ used in the construction of our two-weight set. So with some abuse of notation we have $I_i=(\beta^i)^S$. Note however that this indexing of the $I_i$ may be different from the indexing used earlier in Theorem \ref{main}.
	
	Define the map $\alpha : \Pi \rightarrow \Gamma : \gamma^i \mapsto I_i$. We show that this is an anti-automorphism from $\Pi$ to $\Gamma$. In order to show this consider a hyperplane Tr$(ay)=0$, $a\in\mathbb{F}_{q^n}\setminus\{0\}$, in $\Pi$. Assume that $\gamma^i$ is a point of this hyperplane, that is Tr$(a\gamma^i)=0$. Now the image of this point under $\alpha$ is $I_i=\{ (\beta^i\beta^{j\frac{q^n-1}{q-1}},0) \mid j=0,1,\dots,q^n \}$. We compute that $$\mathrm{Tr}\left(a\left(\beta^i\beta^{j\frac{q^n-1}{q-1}}\right)^{q^n+1}\right)=\mathrm{Tr}\left(a\beta^{(q^n+1)i}\beta^{j\frac{q^{2n}-1}{q-1}}\right)=\beta^{j\frac{q^{2n}-1}{q-1}}\mathrm{Tr}\left(a\gamma^i\right)=0$$
	
	where Tr is the trace map from $\mathbb{F}_{q^n}$ to $\mathbb{F}_q$ and the one but last inequality follows from the fact that $\beta^{j\frac{q^{2n}-1}{q-1}}\in \mathbb{F}_q$. Hence the hyperplane Tr$(ay)=0$ in $\Pi$ gets mapped by $\alpha$ to the elliptic quadric $\mathcal{Q}_a=\mathrm{Tr}(ax^{q^n+1})=0$ in $\Lambda$, that is the point $\mathcal{Q}_a$ of $\Gamma$. It follows that $\alpha$ maps hyperplanes to points in an incidence preserving way and hence that $\alpha$ is an anti-automorphism. We now easily obtain the following description of the two-weight set $\mathcal{C}$:
	
	\begin{theorem}
	Let PG$(3n-1,q)$ be coordinatized as earlier in this section. Then the set $(1,1)^G\cup\Pi$, that is, the orbit of the point $(1,1)$ under the group $G$ together with $\Pi$, is a two-weight set  with weights $\frac{q^{2n}-q^{2n-1}+q^n-q}{q-1}$ and $\frac{q^{2n}-q^{2n-1}-q^{n+1}+2q^n-q}{q-1}$.
	\end{theorem}
	
	{\bf Proof.} This follows from the fact that $\alpha$ as described above is an anti-automorphism and Theorem \ref{main}. \qed
	
\medskip
	
	It follows that the two-weight set $\mathcal{C}$ admits a sharply transitive cyclic group acting on $\mathcal{C}\setminus\Pi$.

	\section{Open Questions and Comments}
	
	We want to end this paper with the following comments and open problems:
	
	\begin{itemize}
	
	\item It is not very hard to see that the two-weight sets constructed through Theorem \ref{main} are not geometric, that is, they are not the union of disjoint $(n-1)$-dimensional subspaces. One can actually see that $\Pi$ is the only $(n-1)$-dimensional subspace contained in $\mathcal{C}$. As a consequence they do not generalize Mathon's construction of a two-weight set of type $(21,30)$ in PG$(5,3)$.
	
	\item Does the two-weight set $\mathcal{C}$ belong to the family of two-weight sets constructed by Bierbrauer and Edel in \cite{Bierbrauer}? This is unclear to the author but would be of interest to know. If they are it would give a nice starting point to find a geometric description of the full Bierbrauer-Edel two-weight sets. As a general open problem it would be of interest to find such geometric description.
	
	\item Is it, under certain conditions, possible to ``glue'' together multiple two-weight sets as constructed in Theorem \ref{main} sharing the same ``base'' $\Pi$ to form larger two-weight sets? This question is inspired by Denniston's construction of maximal arcs where multiple conics sharing the same nucleus are combined to form larger maximal arcs. The author so far did not succeed in making such construction work. 
		
	\item Do our two-weight sets $\mathcal{C}$ in some cases (when $q$ is not a prime) contain smaller two-weight sets? The author so far did not find any sub-two-weight sets. This question is inspired by the fact that every Denniston maximal $q$-arc in PG$(2,q^n)$, $n>1$, $q$ even, contains Denniston maximal $d$-subarcs for every $d\mid q$.

	\end{itemize}

	\section*{Acknowledgement}

This material is based upon work supported by and done while serving at the National Science Foundation. Any opinions, findings, and conclusions or recommendations expressed in this material are those of the author and do not necessarily reflect the views of the National Science Foundation.

\end{document}